\documentclass[12pt]{amsart}
\usepackage[colorlinks=true, pdfstartview=FitV, linkcolor=blue, citecolor=blue, urlcolor=blue]{hyperref}

\usepackage{amssymb,amsmath, amscd}
\usepackage{times, verbatim}
\usepackage{graphicx}
\usepackage[english]{babel}
 \usepackage[usenames, dvipsnames]{color}
\usepackage{amsmath,amssymb,amsfonts}
\usepackage{enumerate}
\usepackage{anysize}
\marginsize{3cm}{3cm}{3cm}{3cm}
\input xy
\xyoption{all}
\usepackage{pb-diagram}
\usepackage[all]{xy}
\usepackage{pgfplots}
\usepackage{tikz}
\usepackage{xcolor}
\input xy
\xyoption{all}

\DeclareFontFamily{OT1}{rsfs}{}
\DeclareFontShape{OT1}{rsfs}{n}{it}{<-> rsfs10}{}
\DeclareMathAlphabet{\mathscr}{OT1}{rsfs}{n}{it}

\begin{document}
\theoremstyle{plain}

\newtheorem{theorem}{Theorem}[section]
\newtheorem{thm}[equation]{Theorem}
\newtheorem{prop}[equation]{Proposition}
\newtheorem{cor}[equation]{Corollary}
\newtheorem{conj}[equation]{Conjecture}
\newtheorem{lemma}[equation]{Lemma}
\newtheorem{definition}[equation]{Definition}
\newtheorem{question}[equation]{Question}

\theoremstyle{definition}
\newtheorem{conjecture}[theorem]{Conjecture}

\newtheorem{example}[equation]{Example}
\numberwithin{equation}{section}

\newtheorem{remark}[equation]{Remark}

\newcommand{\Hecke}{\mathcal{H}}
\newcommand{\Liea}{\mathfrak{a}}
\newcommand{\Cmg}{C_{\mathrm{mg}}}
\newcommand{\Cinftyumg}{C^{\infty}_{\mathrm{umg}}}
\newcommand{\Cfd}{C_{\mathrm{fd}}}
\newcommand{\Cinftyfd}{C^{\infty}_{\mathrm{ufd}}}
\newcommand{\sspace}{\Gamma \backslash G}
\newcommand{\PP}{\mathcal{P}}
\newcommand{\bfP}{\mathbf{P}}
\newcommand{\bfQ}{\mathbf{Q}}
\newcommand{\Siegel}{\mathfrak{S}}
\newcommand{\g}{\mathfrak{g}}
\newcommand{\A}{\mathbb{A}}
\newcommand{\Q}{\mathbb{Q}}
\newcommand{\Gm}{\mathbb{G}_m}
\newcommand{\Nm}{\mathbb{N}m}
\newcommand{\ii}{\mathfrak{i}}
\newcommand{\II}{\mathfrak{I}}

\newcommand{\kk}{\mathfrak{k}}
\newcommand{\nn}{\mathfrak{n}}
\newcommand{\tF}{\tilde{F}}
\newcommand{\p}{\mathfrak{p}}
\newcommand{\m}{\mathfrak{m}}
\newcommand{\bb}{\mathfrak{b}}
\newcommand{\Ad}{{\rm Ad}\,}
\newcommand{\ttt}{\mathfrak{t}}
\newcommand{\frakt}{\mathfrak{t}}
\newcommand{\U}{\mathcal{U}}
\newcommand{\Z}{\mathbb{Z}}
\newcommand{\bfG}{\mathbf{G}}
\newcommand{\bfT}{\mathbf{T}}
\newcommand{\R}{\mathbb{R}}
\newcommand{\ST}{\mathbb{S}}
\newcommand{\h}{\mathfrak{h}}
\newcommand{\bC}{\mathbb{C}}
\newcommand{\C}{\mathbb{C}}

\newcommand{\N}{\mathbb{N}}
\newcommand{\qH}{\mathbb {H}}
\newcommand{\temp}{{\rm temp}}
\newcommand{\Hom}{{\rm Hom}}
\newcommand{\Aut}{{\rm Aut}}
\newcommand{\rk}{{\rm rk}}
\newcommand{\Ext}{{\rm Ext}}
\newcommand{\End}{{\rm End}\,}
\newcommand{\Ind}{{\rm Ind}}
\newcommand{\ind}{{\rm ind}}
\newcommand{\Irr}{{\rm Irr}}
\def\circG{{\,^\circ G}}

\newcommand{\sgn}{\text{Sgn}}
\newcommand{\perm}{\text{Perm}}
\newcommand{\Mod}{\text{mod}}
\newcommand{\Char}{\text{char}}
\newcommand{\F}{\mathscr{F}_{\mu}}
\newcommand{\Diag}{\text {diag }}
\newcommand{\M}{\underline{\mu}}
\newcommand{\T}{\underline{t}}

\def\diag{{\rm diag}}
\def\Ad{{\rm Ad}}
\def\As{{\rm As}}
\def\wG{{\widehat G}}
\def\G{{\rm G}}

\def\SL{{\rm SL}}
\def\PSL{{\rm PSL}}
\def\GSp{{\rm GSp}}
\def\PGSp{{\rm PGSp}}
\def\Sp{{\rm Sp}}
\def\St{{\rm St}}
\def\GU{{\rm GU}}
\def\SU{{\rm SU}}
\def\U{{\rm U}}
\def\GO{{\rm GO}}
\def\GL{{\rm GL}}
\def\PGL{{\rm PGL}}
\def\GSO{{\rm GSO}}
\def\GSpin{{\rm GSpin}}
\def\GSp{{\rm GSp}}
\def\Char{{\Theta}}
\def\Gal{{\rm Gal}}
\def\SO{{\rm SO}}
\def\O{{\rm  O}}
\def\Mod{{\rm  \bmod}}
\def\Sym{{\rm Sym}}
\def\sym{{\rm sym}}
\def\St{{\rm St}}
\def\Sp{{\rm Sp}}
\def\tr{{\rm tr\,}}
\def\ad{{\rm ad\, }}
\def\Ad{{\rm Ad\, }}
\def\rank{{\rm rank\,}}

\subjclass{Primary 11F70; Secondary 22E55}

\title{Character factorizations for representations of ${\rm GL}(n,\C)$ }

\author{Chayan Karmakar}

\address{Indian Institute of Technology Bombay, Powai, Mumbai-400076}
\email{karmakar.baban28@gmail.com}

\maketitle
    {\hfill \today}

\begin{abstract}
We give another proof of a theorem of D. Prasad (Theorem 2, \textit{Israel J. Math.} 2016), which is also a classical result of Littlewood--Richardson (Theorem VI, \textit{Q. J. Math.} 1934). For integers $m,n \ge 2$, this result calculates the character of an irreducible representation of $\GL(mn,\C)$ at diagonal elements with eigenvalues $\omega^{j-1}_nt_i$ for $1 \le i \le m$, $1 \le j \le n$, where $\omega_n=e^{2\pi \imath/n}$, expressing it as a product of certain characters for $\GL(m,\C)$ evaluated at $\underline{t}^n={\rm diag}(t_1^{n},t_{2}^{n},\dots,t_{m}^{n})$. Unlike previous approaches that rely on determinantal identities, our proof utilizes a direct combinatorial cancellation argument within the Weyl group.
\end{abstract}

\tableofcontents

\section{Introduction}
In the work \cite{DP}, D. Prasad established a factorization theorem for characters of $\GL(mn,\C)$ evaluated at specific elements of the diagonal torus. These elements are of the form 
\[
  \T \cdot c_n=\T \otimes c_n :=
  \begin{pmatrix}
    \T & & \\
      & \omega_n\T & \\
    & & \ddots & \\
    & & & \omega_n^{n-1}\T
  \end{pmatrix},
\] 
where $\T={\rm diag}(t_1,t_2,\dots,t_m)$, $\omega_n$ is a primitive $n$-th root of unity, and 
\[
c_n={\rm diag}(1,\omega_n,\dots,\omega^{n-1}_n).
\]
Throughout this paper, we assume that $m,n \geq 2$. 
D. Prasad proved that the character of a finite-dimensional highest weight representation $\pi_{\underline{\lambda}}$ of $\GL(mn,\C)$ with highest weight $\underline{\lambda}$, when evaluated at $\T \cdot c_n$, factorizes into a product of characters of certain highest weight representations of $\GL(m,\C)$ evaluated at $\T^{n}={\rm diag}(t_1^{n},t_2^{n},\dots,t_m^{n})$. 

This result was recently generalized to all classical groups in~\cite{AANK}. It was also observed in~\cite{AANK} that this factorization for type A was originally discovered by D.~E.~Littlewood and A.~R.~Richardson in 1934 (see~\cite{LR}). 
In the paper~\cite{LR}, Littlewood and Richardson calculated the values of S-functions (which coincide with the Weyl characters in type~A) at the roots of the equation $(x^n-1)^m=0$. They showed that these values vanish under conditions identical to those we derive in Proposition~\ref{prop 1}, and in the non-vanishing case, they factorize into a product of S-functions of lower degree.
The approaches in \cite{DP}, \cite{AANK}, and \cite{LR} all rely on direct manipulation of the Weyl character formula expressed as a determinantal identity. The present work aims at giving another proof of their factorization theorem.

We outline our proof-strategy in brief. Let $\underline{\lambda}=(\lambda_1 \geq \lambda_2 \geq \dots \geq \lambda_{mn})$ be a highest weight for $\GL(mn,\C)$ and let $\rho_{mn}$ be half the sum of positive roots. 
First, we observe that for the character to be nonzero at $\T \cdot c_n$, the weight $\underline{\lambda}+\rho_{mn}$ must satisfy a necessary condition: for each $k \in \{0, \dots, n-1\}$, exactly $m$ integers in the weight must be congruent to $k$ modulo $n$ (see Proposition \ref{prop 1}). 
Assuming this condition holds, we replace $\underline{\lambda}+\rho_{mn}$ with a conjugate $\underline{\mu} = \sigma_0 \cdot(\underline{\lambda}+\rho_{mn})$, where $\sigma_0 \in \mathbb{S}_{mn}$ rearranges the entries such that the first $m$ are congruent to $0$ modulo $n$, the next $m$ are congruent to $1$ modulo $n$, and so on. This conjugation affects the Weyl numerator only by a sign.

The core of our argument involves analyzing the Weyl numerator and denominator separately. The Weyl denominator for $\GL(mn,\C)$ at $\T \cdot c_n$ factors (up to a sign) into the Weyl denominator for $\GL(m,\C)^n$ at $\underbrace{(\T^n,\dots,\T^{n})}_{n \ {\rm times}}$, scaled by a constant $E$ and a monomial $(\prod^{m}_{i=1}t_i)^{\frac{n(n-1)}{2}}$. For the Weyl numerator, we sum over the symmetric group $\mathbb{S}_{mn}$ by decomposing it into left cosets of a subgroup $R \cong (\mathbb{S}_{m})^{n}$. We identify a complementary subgroup $C \cong (\mathbb{S}_{n})^{m}$ such that $R \cap C = \{id\}$. We prove that the summation over cosets not having a representative in $CR$ vanishes. Consequently, the sum restricts to cosets having a representative in $CR$, which allows the Weyl numerator to factorize (up to a sign) into a product of characters of a highest weight representation of $\GL(m,\C)^n$ at the diagonal element $\underbrace{(\T^n,\dots,\T^{n})}_{n \ {\rm times}}$, appearing with the same constant $E$ and the same monomial $(\prod^{m}_{i=1}t_i)^{\frac{n(n-1)}{2}}$. Combining the factorized numerator and denominator, these common terms cancel to yield the main theorem (Theorem \ref{4.1} below).

The paper is written in the hope that such manipulations of the Weyl group may be applicable in other contexts, potentially yielding character identities similar to those in \cite{DP}, \cite{LR}, and \cite{AANK}.

\section{Notations and Preliminaries}\label{sec:notation}
We begin with some general notation to be used throughout the paper.
Let $\underline{\zeta}=(\zeta_1,\zeta_2,\ldots,\zeta_{mn}) \in \mathbb{Z}^{mn}$ be an arbitrary weight of the maximal torus $(\mathbb{C}^{\ast})^{mn}$ of $\GL(mn,\C)$. The symmetric group $\mathbb{S}_{mn}$ acts on the weight $\underline{\zeta}$ by permuting its entries. For $\sigma \in \mathbb{S}_{mn}$, we define:
\[
\sigma \cdot \underline{\zeta}=(\zeta_{\sigma^{-1}(1)},\zeta_{\sigma^{-1}(2)},\ldots,\zeta_{\sigma^{-1}(mn)}).
\]

Define the matrix map $M$ by:
\begin{equation}\label{defn 1}
M(\underline{\zeta}):=(\zeta_{i,j})_{1 \leq i \leq n, 1 \leq j \leq m} \in \mathbb{Z}^{n \times m}, \quad \text{where} \quad \zeta_{i,j}:=\zeta_{(i-1)m+j}.
\end{equation}

Also define the action of $\mathbb{S}_{mn}$ on the matrix entries via:
\begin{equation}\label{defn 2}
\zeta^{\sigma}_{i,j}:=\zeta_{\sigma^{-1}((i-1)m+j)}.
\end{equation}
Then the map $M$ is $\mathbb{S}_{mn}$-equivariant: $M(\sigma \cdot \underline{\zeta}) = (\zeta^{\sigma}_{i,j}) =: M(\underline{\zeta})^{\sigma}$.

\smallskip

Finally, given a permutation $\sigma \in \mathbb{S}_{mn}$ and a weight $\underline{\zeta}$, we define the character $\sigma \cdot \underline{\zeta}$ on an arbitrary element $x={\rm diag}(x_1, \dots, x_{mn}) \in (\mathbb{C}^*)^{mn}$ by:
\begin{equation}\label{defn 3}
(\sigma \cdot \underline{\zeta})(x) := \prod_{k=1}^{mn} x_k^{\zeta_{\sigma^{-1}(k)}}.
\end{equation}
In particular, when evaluated at the specific torus element $\T \cdot c_n$ (where the diagonal entries are given by $x_{(i-1)m+j} = \omega_n^{i-1} t_j$), this becomes:
\[
(\sigma \cdot \underline{\zeta}) (\T \cdot c_n)=\prod_{i=1}^n \prod_{j=1}^m (\omega^{i-1}_n t_j)^{\zeta^{\sigma}_{i,j}}.
\] 
\section{Calculation of the Weyl Denominator $A_{\rho}$}

The Weyl character formula for the highest weight representation $S_{(\underline{\lambda})}\C^{n}$ of $\GL(n,\C)$ with highest weight
$\underline{\lambda}=(\lambda_1 \geq \lambda_2 \geq \cdots \geq \lambda_n)$ is given by 
\begin{equation}\label{defn 4}
\Theta_{\underline{\lambda}}=\frac{\sum_{\sigma \in \mathbb{S}_{n}} (-1)^{\sigma}\sigma \cdot(\underline{\lambda}+\rho_{n})}{\sum_{\sigma \in \mathbb{S}_{n}}(-1)^{\sigma} \sigma \cdot\rho_{n}}=\frac{A_{\underline{\lambda}+\rho_n}}{A_{\rho_n}}, 
\end{equation}
where the Weyl numerator is given by \( A_{\underline{\lambda}+\rho_n}=\sum_{\sigma \in \mathbb{S}_{n}} (-1)^{\sigma} \sigma \cdot(\underline{\lambda}+\rho_{n}), \) and the Weyl denominator is given by \( A_{\rho_n}=\sum_{\sigma \in \mathbb{S}_{n}} (-1)^{\sigma}\sigma \cdot\rho_n.\) Note that the Weyl denominator $A_{\rho_n}$ at $\underline{x}={\rm diag}(x_1,x_2,\cdots,x_n) \in (\C^{\ast})^{n}$ yields the Vandermonde determinant \[ \prod_{i<j}(x_i-x_j). \]

Our work involves factorizing both the Weyl numerator and the denominator. While the Weyl denominator is a special case of the Weyl numerator (for $\underline{\lambda}=0$) \textendash \ so that it suffices to factorize the latter \textendash \ the Weyl denominator $A_{\rho_{mn}}(\T \cdot c_n)$ is a much simpler expression, and its factorization is easy enough, so we begin by factoring it in this section.   
\begin{align*}
A_{\rho_{mn}}(\T \cdot c_n,x)& \stackrel{(1)}{=} \prod_{ \{ k < l ,\hspace{0.001cm} i \geq j \}}(\omega_n^{k}t_i-\omega_n^{l}t_j)  \times  \prod_{ \{   i < j, \hspace{0.01cm} k \leq l \}}(\omega_n^{k}t_i-\omega_n^{l}t_j) \times , \\
& \stackrel{(2)}{=}\prod_{k < l,\hspace{0.01cm}}\prod^{m}_{i=1}(\omega_n^{k}-\omega_n^{l})t_i  \ \ \ \ \ \times \prod_{ \{ k < l ,\hspace{0.001cm} i > j \}}(\omega_n^{k}t_i-\omega_n^{l}t_j)  \times \prod_{ \{   i < j, \hspace{0.01cm} k \leq l \}}(\omega_n^{k}t_i-\omega_n^{l}t_j), \\
&= \hspace{1.02cm} A  \hspace{2.20cm}  \times  \ \ B \hspace{3.20cm} \times C. 
\end{align*} 

We have
 \begin{align}
 A&=\prod_{k<l} (\omega_n^{k}-\omega_n^{l})^{m} \cdot (\prod_{i=1}^{m}t_i)^{\frac{n(n-1)}{2}}.
 \end{align}
 Now we will evaluate $B \times C$. We have
 $$
\begin{aligned}
B \times C & =\prod_{i<j, k \leq l}\left(\omega_{n}^{k} t_{i}-\omega_{n}^{l} t_{j}\right) \cdot \prod_{j<i, k<l}\left(\omega_{n}^{k} t_{i}-\omega_{n}^{l} t_{j}\right) \\
& =\prod_{i<j, k \leq l}\left(\omega_{n}^{k} t_{i}-\omega_{n}^{l} t_{j}\right) \cdot \prod_{1 \leq i<j \leq m, 0 \leq l<k \leq n-1}\left(\omega_{n}^{l} t_{j}-\omega_{n}^{k} t_{i}\right) \\
& =\prod_{i<j} \prod_{k, l=0}^{n-1}\left(\omega_{n}^{k} t_{i}-\omega_{n}^{l} t_{j}\right) \cdot(-1)^{\binom{m}{2}\binom{n}{2}} \\
& =(-1)^{\binom{m}{2}\binom{n}{2}} \prod_{i<j} \prod_{k=0}^{n-1}\left(\omega_{n}^{k n} \prod_{l^{\prime}=0}^{n-1}\left(t_{i}-\omega_{n}^{l^{\prime}} t_{j}\right)\right) \\
& =(-1)^{\binom{m}{2}\binom{n}{2}} \prod_{i<j} \prod_{k=0}^{n-1}\left(t_{i}^{n}-t_{j}^{n}\right) 
\end{aligned}
$$
Therefore 
\begin{align}
\nonumber
A_{\rho_{mn}}(\T \cdot c_n)&=A \times B \times C, \\
&=(-1)^{\frac{mn(m-1)(n-1)}{4}} \times \prod_{k<l} (\omega_n^{k}-\omega_n^{l})^{m} \times \prod_{i < j}(t_i^{n}-t_j^{n})^{n} \times (\prod_{i=1}^{m}t_i)^{\frac{n(n-1)}{2}}.\label{4}
\end{align}

\section{Calculation of the Weyl numerator}
The proof of the main theorem of this paper will go through several lemmas and two propositions. 
We first start with a counting lemma about the symmetric group which may be of independent interest. 

\begin{lemma}\label{lem 1}
 For the weight $\underline{\zeta}=(\zeta_1,\zeta_2,\cdots,\zeta_{mn})$ of $\GL(mn,\C)$, let us consider $M(\underline{\zeta})=(\zeta_{i,j})_{ij}$, where $\zeta_{i,j}=\zeta_{(i-1)m+j}$, defined in Section ~\ref{sec:notation}. Let $R \cong (\mathbb{S}_{m})^{n}$ be the subgroup of $\mathbb{S}_{mn}$ which fixes the rows of $M(\underline{\zeta})$ and $C \cong (\mathbb{S}_{n})^{m}$ be the subgroup of $\mathbb{S}_{mn}$ which fixes the columns (here $R$ permutes the elements of each row of $M(\underline{\zeta})$ and similarly $C$ permutes the elements of each column). We clearly have $R \cap C=\{ id_{mn}\}.$ Now a permutation $\tau  \in CR$ if and only if $\tau$ takes no two entries of the same row to the same column.
\end{lemma} 
\begin{proof}
We may write
\[
M(\underline{\zeta}) = (\zeta_{i,j}) =
\begin{pmatrix}
\zeta_{1} & \zeta_{2} & \zeta_{3} & \dots & \zeta_{m-1} & \zeta_{m} \\
\zeta_{m+1} & \zeta_{m+2} & \zeta_{m+3} & \dots & \zeta_{2m-1} & \zeta_{2m} \\
\vdots & \vdots & \vdots & \ddots & \vdots & \vdots \\
\zeta_{(n-1)m+1} & \zeta_{(n-1)m+2} & \zeta_{(n-1)m+3} & \dots & \zeta_{nm-1} & \zeta_{mn}
\end{pmatrix}.
\]

If $\tau \in CR$, then by definition $\tau$ does not send two entries from the 
same row to entries in the same column.  We count the number of permutations 
with this property.

For the first row, the $m$ entries must be sent to $m$ distinct columns.  The 
number of ways to do this is
\[
mn \cdot (mn-n) \cdots (mn-(m-1)n) = m! \, n^{m}.
\]
After these entries are placed, there remain $mn-m$ entries.  The number of ways 
to place the entries of the second row into distinct columns is $m! \, (n-1)^{m}$.  
Proceeding similarly, the number of ways to place the entries of the $k$-th row 
into distinct columns is
\[
m! \, (n-k+1)^{m}.
\]

Therefore, the total number of permutations that do not send two entries of the 
same row to the same column is
\[
\prod_{k=1}^{n} m! \, (n-k+1)^{m}
= (m!)^{n} \, (n!)^{m}
= \#CR.
\]
Since every element of $CR$ satisfies the stated property, and the number of such 
permutations equals $\#CR$, these are exactly the elements of $CR$.
\end{proof}

\noindent
Note that the above lemma depends only on the combinatorial structure of the 
matrix $M(\underline{\zeta})$ and the actions of the subgroups $R$ and $C$, and is 
therefore independent of the particular choice of the weight $\underline{\zeta}$.
For any $w \in \mathbb{S}_{mn}$, we have 
\[
\sum_{\sigma \in \mathbb{S}_{mn}}(-1)^{\sigma }\sigma w \cdot (\underline{\lambda}+\rho_{mn})=(-1)^{w}A_{\underline{\lambda}+\rho_{mn}}.
\]
So we can always work with a suitable conjugate of $\underline{\lambda}+\rho_{mn}$.  
We will take the Weyl numerator $A_{\underline{\lambda}+\rho_{mn}}$ to be $\sum_{\sigma \in \mathbb{S}_{mn}}(-1)^{\sigma}\sigma \cdot\underline{\mu}$, where $\underline{\mu}$ is a conjugate of $\underline{\lambda}+\rho_{mn}$ by an element of the Weyl group $\mathbb{S}_{mn}$ so that the entries of $\underline{\lambda}+\rho_{mn}$ which are congruent to $0$ modulo $n$ come first in decreasing order, then come those integers which are congruent to $1$ modulo $n$ in decreasing order, and so on.
 
We prove the necessary condition on $\underline{\mu}=(\mu_1,\mu_2,\dots,\mu_{mn})$ for the Weyl numerator $A_{\underline{\mu}}$ to be nonzero in the next proposition.

\begin{prop}\label{prop 1}
If the entries of $\underline{\mu}$ do not satisfy the condition that exactly $m$ integers represent each residue class modulo $n$, then the Weyl numerator $A_{\underline{\mu}}(\T \cdot c_n)$ is zero.
\end{prop}

\begin{proof}
Recall the notations from Section \ref{sec:notation}. The symmetric group $\mathbb{S}_{mn}$ acts on $M(\underline{\mu}) = (\mu_{i,j})_{i,j}$ as defined in equation \eqref{defn 1}.

Under the assumption on $\underline{\mu}$ in the proposition, the entries are not equidistributed among the residue classes modulo $n$. By the Pigeonhole Principle, for every $\sigma \in \mathbb{S}_{mn}$, there must be at least one column—say, the $k_{\sigma}$-th column of $\sigma \cdot M(\underline{\mu})$—that contains two entries $\mu^{\sigma}_{a,k_{\sigma}}$ and $\mu^{\sigma}_{b,k_{\sigma}}$ (assume without loss of generality that $a < b$) satisfying
\[
\mu^{\sigma}_{a,k_{\sigma}} \equiv \mu^{\sigma}_{b,k_{\sigma}} \pmod n.
\]

Consider the transposition
\[
p_{\sigma}
=\bigl(\sigma^{-1}((a-1)m+k_{\sigma})\ \ 
\sigma^{-1}((b-1)m+k_{\sigma})\bigr).
\]
For any $\sigma \in \mathbb{S}_{mn}$, and for indices $1 \leq j \leq m$ and $0 \leq r_1 < r_2 \leq n-1$, we define the set:
\begin{equation}
P_{\sigma}
=\Bigl\{\, \bigl(\sigma^{-1}(r_1m+j)\;\; \sigma^{-1}(r_2m+j)\bigr) \in \mathbb{S}_{mn} \ \Big|\ 
\mu_{\sigma^{-1}(r_1m+j)} \equiv \mu_{\sigma^{-1}(r_2m+j)} \pmod n \Bigr\}.
\end{equation}
For each odd permutation $\sigma$, define
$G_{\sigma}=\langle P_{\sigma}\rangle$, the subgroup of $\mathbb{S}_{mn}$ generated by $P_{\sigma}$.
Note that since $p_{\sigma} \in P_{\sigma}$ (by setting $j = k_{\sigma}$, $r_1 = a-1$, and $r_2 = b-1$), the set $P_{\sigma}$ is nonempty.

Any transposition $t_{\sigma} \in P_{\sigma}$ permutes only those entries of $\sigma \cdot M(\underline{\mu})$ that lie in the same residue class and occur in the same column. As a consequence, we have 
\begin{equation}\label{eq 4.2}
G_{\sigma}=G_{\sigma t_{\sigma}},
\end{equation}
for each $t_{\sigma} \in P_{\sigma}$. 

Since each transposition in $P_{\sigma}$ swaps two entries in the same column 
of $\sigma\cdot M(\underline{\mu})$ that are congruent modulo $n$, it follows that 
\begin{equation}\label{eq 4.3}
\sigma \tau\cdot \underline{\mu}(\T \cdot c_n)=\sigma \cdot \underline{\mu}(\T \cdot c_n), \quad \forall \tau \in G_{\sigma}.
\end{equation}
From equation \eqref{eq 4.3}, considering the sum over the left coset $\sigma G_{\sigma}$, we deduce:
\begin{align}
\nonumber
\sum_{\tau \in \sigma G_{\sigma}}(-1)^{\tau}\tau \cdot \underline{\mu}(\T \cdot c_n)&=\sum_{\rho \in G_{\sigma}}(-1)^{\sigma\rho}\sigma\rho \cdot \underline{\mu}(\T \cdot c_n) \\
\nonumber
&=\sum_{\rho \in G_{\sigma}}(-1)^{\sigma}(-1)^{\rho}\sigma\cdot \underline{\mu}(\T \cdot c_n) \\
\nonumber
&=(-1)^{\sigma}\sigma\cdot \underline{\mu}(\T \cdot c_n) \sum_{\rho \in G_{\sigma}}(-1)^{\rho} \\
&=0.\label{eq 4.5}
\end{align}
The last equality holds because $G_{\sigma}$ contains transpositions (which are odd), so half the elements are even and half are odd.

Next, we show that the collection $\{\sigma G_{\sigma} \mid \sigma \text{ is an odd permutation}\}$ covers the entire symmetric group $\mathbb{S}_{mn}$.
If $\sigma$ is an odd permutation, clearly $\sigma \in \sigma G_{\sigma}$. Let $\sigma$ be an even permutation. Since $P_{\sigma}$ is nonempty, there exists a transposition $t_{\sigma} \in P_{\sigma}$. We can write $\sigma=(\sigma t_{\sigma}) t_{\sigma}$. Note that $\sigma t_{\sigma}$ is an odd permutation.
From \eqref{eq 4.2}, $G_{\sigma t_{\sigma}} = G_{\sigma}$. Thus,
\[ \sigma = (\sigma t_{\sigma}) t_{\sigma} \in (\sigma t_{\sigma}) G_{\sigma} = (\sigma t_{\sigma}) G_{\sigma t_{\sigma}}. \]
Therefore, every even permutation belongs to the coset of an odd permutation.
We obtain
\[
\bigcup_{\sigma \text{ odd}} \sigma G_{\sigma}
= \mathbb{S}_{mn}.
\]   
Note that the above union is not disjoint. We will show that if $\sigma G_{\sigma} \cap \tau G_{\tau} \neq \emptyset$, then we have $\sigma G_{\sigma}=\tau G_{\tau}$. Suppose that $\sigma G_{\sigma} \cap \tau G_{\tau} \neq \emptyset$. Then we have $\sigma u_{\sigma}=\tau u_{\tau}$, where $u_{\sigma} \in G_{\sigma}$ and $u_{\tau} \in G_{\tau}$. Now from \eqref{eq 4.2}, it follows that $G_{\sigma}=G_{\sigma u_{\sigma}}=G_{\tau u_{\tau}}=G_{\tau}$. Therefore $\sigma \tau^{-1} \in G_{\sigma}=G_{\tau}$. Hence we obtain 
$\sigma G_{\sigma}=\tau G_{\tau}$.

To handle the overlap, let us define an equivalence relation on the collection of cosets generated by odd permutations: $\sigma G_{\sigma} \sim \tau G_{\tau}$ if $\sigma G_{\sigma}=\tau G_{\tau}$.
Let $\{\sigma_1 G_{\sigma_1}, \dots, \sigma_k G_{\sigma_k}\}$ be a set of disjoint cosets whose union is $\mathbb{S}_{mn}$.
Finally, from \eqref{eq 4.5} we obtain:
\begin{align}
\nonumber
\sum_{\nu \in \mathbb{S}_{mn}}
(-1)^{\nu}\,\nu \cdot\underline{\mu}(\T \cdot c_n)
&=\sum^{k}_{i=1}  
\bigg(\sum_{\nu_i \in \sigma_i G_{\sigma_i}}(-1)^{\nu_i}\nu_i \cdot\underline{\mu}(\T \cdot c_n) \bigg) \\
&=0.
\label{eqn 4.5}
\end{align}
\end{proof}

\begin{remark}
Proposition~\ref{prop 1} is analogous in spirit to Theorem~VI of \cite{LR}. 
In \cite{LR}, the authors prove the vanishing of the S-function (which corresponds to the Weyl character in type A) by expressing the character formula as a determinant and applying a Laplace expansion to show that the determinant is zero when the residue condition is not met.

In contrast, our proof of Proposition~\ref{prop 1} relies on a direct combinatorial cancellation argument within the Weyl group summation. 
Specifically, if $\underline{\mu}$ does not satisfy the condition that exactly $m$ entries represent each residue class modulo $n$, the Pigeonhole Principle guarantees that for every $\sigma \in \mathbb{S}_{mn}$, there is at least one column of $\sigma \cdot M(\M)$ containing two entries with the same residue modulo $n$.
We then construct for every $\sigma$ a subgroup $G_{\sigma}$ generated by transpositions that swap entries in the same column of $\sigma \cdot M(\underline{\mu})$ having the same residue modulo $n$. 
Since these transpositions are odd permutations that stabilize the term $\sigma \cdot \underline{\mu}(\T \cdot c_n)$, the summation over the coset $\sigma G_{\sigma}$ vanishes.

Crucially, we show that the collection of these cosets $\{\sigma G_{\sigma}\}$, where $\sigma$ varies over odd permutations, covers the entire symmetric group $\mathbb{S}_{mn}$. Although these cosets are not disjoint, we define an equivalence relation to partition the group into disjoint unions of such cosets, thereby ensuring that the entire Weyl numerator vanishes.
This establishes the result via combinatorial necessity rather than determinantal identities.
\end{remark}    
Given Proposition \ref{prop 1}, and since we will be proving our main result up to a sign $\pm 1$, from now on we will take the Weyl numerator $A_{\underline{\lambda}+\rho_{mn}}$ to be $\sum_{w \in \mathbb{S}_{mn}}(-1)^{w} w \cdot\M$, where $\M$ is a conjugate of $\underline{\lambda}+\rho_{mn}$ by an element of the Weyl group $\mathbb{S}_{mn}$ so that it satisfies the property: 
\begin{equation}\label{4.9}
 P\#=\small
 \left\{\begin{array}{lr}
        \text {the first $m$ entries of $\M$ are congruent to $0$ modulo $n$}, \\
        \text{ the next $m$ entries are congruent to $1$ modulo $n$ etc.} 
        \end{array}\right\}.     
\end{equation} 
The next lemma will determine the conditions on the left cosets $\tau R$ for which the summation $\sum_{\sigma \in R}(-1)^{\tau \sigma} (\tau \sigma)\cdot \M(\T \cdot c_n)$ is zero. 

\begin{lemma}\label{lem 2}
If $\underline{\mu}$ satisfies property $P\#$, then with the notation as in Lemma~\ref{lem 1},
\begin{equation}
\sum_{\sigma \in R}(-1)^{\tau \sigma}\tau \sigma \cdot\underline{\mu}(\T \cdot c_n)=0, \quad \text{for all } \tau \not \in CR. 
\end{equation}
\end{lemma}

\begin{proof}
As in Proposition~\ref{prop 1}, we write $\underline{\mu}$ in the matrix form $M(\underline{\mu})=(\mu_{i,j})_{i,j}$, where $\mu_{i,j}=\mu_{(i-1)m+j}$. Thus, $M(\underline{\mu})$ is an $n \times m$ matrix with the subgroup $C$ permuting the entries along the columns and the subgroup $R$ permuting the entries along the rows (see Lemma~\ref{lem 1}). Recall that $\mathbb{S}_{mn}$ acts on $M(\underline{\mu})=(\mu_{i,j})_{i,j}$ as in \eqref{defn 1}.

From Lemma~\ref{lem 1}, if $\tau \notin CR$, then $\tau$ sends at least two entries from some row of $M(\underline{\mu})$—say, the $a$-th row—to the same $b$-th column of $\tau\cdot M(\underline{\mu})$. Therefore, there exist entries $\mu^{\tau}_{k_1,b}$ and $\mu^{\tau}_{k_2,b}$ in the matrix $\tau \cdot M(\underline{\mu})$ such that 
\begin{equation}
 \mu^{\tau}_{k_1,b} \equiv \mu^{\tau}_{k_2,b} \equiv (a-1) \pmod n.
\end{equation}
This congruence holds because $\underline{\mu}$ satisfies property $P\#$, meaning all entries originating from the $a$-th row share the residue $(a-1)$ modulo $n$. Because the entries $\mu^{\tau}_{k_1,b}$ and $\mu^{\tau}_{k_2,b}$ actually come from the $a$-th row, applying any $\sigma \in R$ will only change their initial positions within that same $a$-th row. Thus, we have:
\begin{equation}\label{Cong. relation}
 \mu^{\tau \sigma}_{k_1,b} \equiv \mu^{\tau \sigma}_{k_2,b} \pmod n. 
\end{equation}

Consider the transposition
\begin{equation}
\sigma_0 = \bigl( \tau^{-1}((k_1-1)m+b) \;\; \tau^{-1}((k_2-1)m+b) \bigr).
\end{equation}
Since $\sigma_0$ swaps two elements that both originally come from the $a$-th row of $M(\underline{\mu})$, it is a row-preserving permutation, and hence $\sigma_0 \in R$. For any $\sigma \in R$, applying $\sigma_0$ yields:
\begin{equation}\label{eqn 4.9}
 \mu^{\tau \sigma_0\sigma}_{i,j}=
 \begin{cases}
\mu^{\tau\sigma}_{k_2,b}, & \text{if } i=k_1, j=b, \\ 
\mu^{\tau\sigma}_{k_1,b}, & \text{if } i=k_2, j=b, \\
\mu^{\tau\sigma}_{i,j}, & \text{otherwise}.
 \end{cases}
\end{equation} 

Now we calculate $\tau \sigma_0\sigma \cdot \underline{\mu}(\T \cdot c_n)$. Using equations \eqref{Cong. relation} and \eqref{eqn 4.9}, we separate the terms in the $b$-th column and deduce:
\begin{align}
\nonumber
\tau\sigma_{0}\sigma \cdot \underline{\mu}(\T \cdot c_n) &= \prod^{m}_{j=1}\left(\prod^{n}_{i=1}(\omega_n^{i-1}t_j)^{\mu^{\tau\sigma_{0}\sigma}_{i,j}}\right) \\
\nonumber
&= \prod_{j \neq b}\left(\prod^{n}_{i=1}(\omega_n^{i-1}t_j)^{\mu^{\tau\sigma}_{i,j}}\right) \times \left(\prod^{n}_{i=1}(\omega_n^{i-1}t_b)^{\mu^{\tau\sigma_{0}\sigma}_{i,b}}\right) \\
\nonumber
&= \prod_{j \neq b}\left(\prod^{n}_{i=1}(\omega_n^{i-1}t_j)^{\mu^{\tau\sigma}_{i,j}}\right) \times \left(\prod_{i \neq k_1,k_2}(\omega_n^{i-1}t_b)^{\mu^{\tau\sigma}_{i,b}} \right) \\
\nonumber
&\qquad \times (\omega_n^{k_1-1}t_b)^{\mu^{\tau\sigma}_{k_2,b}} \times (\omega_n^{k_2-1}t_b)^{\mu^{\tau\sigma}_{k_1,b}}. 
\end{align}
Because $\mu^{\tau \sigma}_{k_1,b} \equiv \mu^{\tau \sigma}_{k_2,b} \pmod n$, their difference is a multiple of $n$. Since $\omega_n^n = 1$, the value of the last product is invariant under swapping these two exponents. Thus, the last two factors can be rewritten as $(\omega_n^{k_1-1}t_b)^{\mu^{\tau\sigma}_{k_1,b}} \times (\omega_n^{k_2-1}t_b)^{\mu^{\tau\sigma}_{k_2,b}}$, which gives:
\begin{align}
\tau\sigma_{0}\sigma \cdot \underline{\mu}(\T \cdot c_n) &= \tau\sigma \cdot\underline{\mu}(\T \cdot c_n). \label{eqn 4.11}
\end{align}

From equation \eqref{eqn 4.11}, we observe that substituting $\sigma$ with $\sigma_0 \sigma$ in the sum yields:
\begin{align*}
\sum_{\sigma \in R}(-1)^{\tau \sigma}\tau \sigma\cdot \underline{\mu}(\T \cdot c_n) &= \sum_{\sigma \in R}(-1)^{\tau \sigma_0 \sigma}\tau \sigma_0 \sigma \cdot\underline{\mu}(\T \cdot c_n) \\
&= (-1)^{\sigma_0}\sum_{\sigma \in R}(-1)^{\tau \sigma} \tau \sigma\cdot\underline{\mu}(\T \cdot c_n) \\
&= -\sum_{\sigma \in R}(-1)^{\tau \sigma}\tau \sigma \cdot\underline{\mu}(\T \cdot c_n).
\end{align*}
Therefore, $\sum_{\sigma \in R}(-1)^{\tau \sigma}\tau \sigma \cdot\underline{\mu}(\T \cdot c_n)=0$ for all $\tau \not \in CR$.
\end{proof}
Next we will calculate $ \sum_{\sigma \in R}(-1)^{\tau \sigma} \tau \sigma \cdot\M(\T \cdot c_n)$, where $\tau \in CR$, in the following  lemma and a proposition.

\begin{lemma}\label{lem 3}
Let $\M$ be a character of the maximal torus $(\C^{\ast})^{mn}$ of $\GL(mn,\C)$ satisfying property $P\#$. For any permutation $g \in C$ (see Lemma~\ref{lem 1}), there exists an $n$-th root of unity $C_{g}(\M) \in \C$ such that:
\begin{itemize}
    \item [(i)] $g \cdot\M(\T \cdot c_n) = \omega^{-m\big(\sum^{n-1}_{k=0}k^{2}\big)}_{n} \cdot C_{g}(\M) \cdot \M(\T \cdot c_n).$
    
    \item [(ii)] For any $g \in C$ and $\sigma \in R$, we have $C_{g\sigma}(\M) = C_{g}(\M)$, and consequently: 
    \[
    g\sigma \cdot \M(\T \cdot c_n) = \omega^{-m\big(\sum^{n-1}_{k=0}k^{2}\big)}_{n} \cdot C_{g}(\M) \sigma \cdot\M(\T \cdot c_n).
    \]
    
    \item [(iii)] Let $C=\prod^{m}_{j=1}P_{j}$, where $P_{j} \cong \mathbb{S}_n$ is the group of permutations of the $j$-th column of $M(\M)=(\mu_{i,j})$. For any $g_j \in P_{j}$, we have:
    \[
    g_{j}\cdot\M^{j}(t_j,\omega_nt_j,\dots,\omega^{n-1}_nt_j) = C_{g_{j}}(\M^{j}) \, t^{\sum^{n}_{i=1}\mu^{g_{j}}_{i,j}}_{j},
    \] 
    where $\M^{j}=(\mu_{1,j}, \dots, \mu_{n,j})^T$ denotes the $j$-th column of $M(\M)$, viewed as a character of $\GL(n,\C)$, and the constant $C_{g_{j}}(\M^{j}) \in \C$ is an $n$-th root of unity.
    
    \item[(iv)] If $g=\prod^{m}_{j=1}g_{j}$ with $g_{j} \in P_{j}$, then $C_{g}(\M) = \prod^{m}_{j=1}C_{g_{j}}(\M^{j})$.
    
    \item[(v)] The alternating sum factors over the columns:
    \[
    \sum_{g \in C}(-1)^{g}C_{g}(\M) = \prod^{m}_{j=1}\left[ \sum_{g_{j} \in P_{j}}(-1)^{g_j}C_{g_{j}}(\M^{j}) \right].
    \]
    
    \item[(vi)] For the permutation group $P_{j}$ of the $j$-th column, the sum yields a Vandermonde determinant:
    \begin{equation*}
     \sum_{g_{j} \in P_{j}}(-1)^{g_{j}}C_{g_{j}}(\M^{j}) = \pm \det
     \begin{pmatrix}
     1 & 1 & \dots & 1 \\
     1 & \omega_{n} & \dots & \omega^{n-1}_{n} \\
     1 & \omega^{2}_{n} & \dots & \omega^{2(n-1)}_{n} \\
     \vdots & \vdots & \ddots & \vdots \\
     1 & \omega^{n-1}_{n} & \dots & \omega^{(n-1)^{2}}_{n} 
     \end{pmatrix}
     = \pm \prod_{0 \le k < l \le n-1}(\omega^{l}_{n}-\omega^{k}_{n}). 
    \end{equation*}
\end{itemize}
\end{lemma}

\begin{proof}
The elements of the subgroup $C$ permute the entries of each column of $M(\M)$. Using this, we deduce:
\begin{align}
g \cdot \M(\T \cdot c_n) &= \prod^{m}_{j=1}\left( \prod^{n}_{i=1}(\omega^{i-1}_nt_j)^{\mu^{g}_{i,j}}\right) \nonumber \\
&= \omega^{\sum^{m}_{j=1}\sum^{n}_{i=1}(i-1)\mu^{g}_{i,j}}_{n} \times \prod^{m}_{j=1}t^{\sum^{n}_{i=1}\mu^{g}_{i,j}}_j \nonumber \\
&= \omega^{\sum^{m}_{j=1}\sum^{n}_{i=1}(i-1)\mu^{g}_{i,j}}_{n} \times \prod^{m}_{j=1}t^{\sum^{n}_{i=1}\mu_{i,j}}_j. \label{eqn 4.13}    
\end{align}

Similarly, evaluating the identity permutation yields:
\begin{align}
\M(\T \cdot c_n) &= \prod^{m}_{j=1}\left( \prod^{n}_{i=1}(\omega^{i-1}_nt_j)^{\mu_{i,j}}\right) \nonumber \\
&= \omega^{\sum^{m}_{j=1}\sum^{n}_{i=1}(i-1)\mu_{i,j}}_{n} \times \prod^{m}_{j=1}t^{\sum^{n}_{i=1}\mu_{i,j}}_j \nonumber \\
&= \omega^{m \sum^{n}_{i=1}(i-1)^{2}}_{n} \times \prod^{m}_{j=1}t^{\sum^{n}_{i=1}\mu_{i,j}}_j. \label{eq 4.14}
\end{align}

Part (i) follows directly by dividing equation \eqref{eqn 4.13} by equation \eqref{eq 4.14}. The products involving the parameters $t_j$ are identical in both expressions and cancel completely. By factoring out the permutation-independent scalar $\omega_{n}^{-m \sum_{k=0}^{n-1}k^{2}}$, we isolate the effect of the permutation $g$ and define the remaining constant as:
\begin{equation}\label{defn Const 1}
C_{g}(\M) = \omega^{\sum^{m}_{j=1} \sum^{n}_{i=1}(i-1)\mu^{g}_{i,j}}_n.
\end{equation}

To prove part (ii), observe that since the subgroup $R$ permutes entries strictly along the rows of $M(\M)$, applying $\sigma \in R$ does not change the row from which any element originated. Thus, the index $i$ in the exponent remains unchanged:
\[
\sum^{m}_{j=1} \sum^{n}_{i=1}(i-1)\mu^{g\sigma}_{i,j} = \sum^{m}_{j=1} \sum^{n}_{i=1}(i-1)\mu_{\sigma^{-1}g^{-1}((i-1)m+j)} \equiv \sum^{m}_{j=1} \sum^{n}_{i=1}(i-1)\mu^{g}_{i,j} \pmod n.
\]
Therefore, $C_{g\sigma}(\M) = C_{g}(\M)$ for all $\sigma \in R$. The transformation equation follows directly from (i).

For part (iii), viewing the $j$-th column $\M^{j}=(\mu_{1,j}, \dots, \mu_{n,j})^T$ as a character of $\GL(n,\C)$ evaluated at the torus element $\text{diag}(t_j,\omega_nt_j,\dots,\omega^{n-1}_{n}t_j)$, we have:
\[
g_j \cdot \M^{j}(t_j,\dots,\omega^{n-1}_{n}t_j) = C_{g_{j}}(\M^{j}) \, t^{\sum^{n}_{i=1}\mu_{i,j}}_{j},
\]
where we define
\begin{equation}\label{defn Const 2}
C_{g_{j}}(\M^{j}) = \prod^{n}_{i=1}(\omega^{i-1}_n)^{\mu^{g_{j}}_{i,j}}.
\end{equation}

Part (iv) follows immediately by comparing equation \eqref{defn Const 1} with the product of the terms defined in equation \eqref{defn Const 2}. Part (v) is a direct consequence of (iv) and the direct product structure $C=\prod^{m}_{j=1}P_{j}$. 

Finally, for part (vi), because $\M$ satisfies property $P\#$, the $n$ entries of the $j$-th column $\M^{j}$ are mutually incongruent modulo $n$. Consequently, as $g_j$ varies over the symmetric group $P_j$, the alternating sum $\sum_{g_j \in P_j}(-1)^{g_j}C_{g_j}(\M^{j})$ expands, up to a sign, as the determinant of the standard Vandermonde matrix generated by the $n$-th roots of unity $1, \omega_n, \dots, \omega_n^{n-1}$. This overall sign $\pm 1$ depends on the initial ordering of the residues in the column. This completes the proof.
\end{proof}

\begin{prop}\label{prop 3.6}
For $\underline{\mu}$ a character of $(\C^{\ast})^{mn}$ with the property $P\#$ , and $g \in C$,
\[
\sum_{\sigma \in R}(-1)^{g  \sigma} g \sigma \cdot \M(\T \cdot c_n)=(-1)^{g}C_{g}(\M) \cdot \left (\prod_{i=1}^{n} S_i(\T^{n})\right) \cdot (\prod_{j=1}^{m}t_j)^{\frac{n(n-1)}{2}},
\]
where $S_{i}$ is the Weyl numerator for $\GL(m,\C)$ with highest weight $\eta_{i}=\left ([\underline{\mu_{i}}-i+1]/n-\rho_{m}\right),$  $\underline{\mu_{i}}$ is the $i$-th row of $M(\M)=(\mu_{i,j})$.
\end{prop}
\begin{proof}
By Lemma \ref{lem 3}, it suffices to prove the above identity for $g=1$, in which case the group $R$ is the product $\underbrace{\mathbb{S}_{m} \times \mathbb{S}_{m} \times \cdots \times \mathbb{S}_{m}}_{n \ {\rm times}}$, where the $i$-th copy of $\mathbb{S}_{m}$ $(1 \leq i \leq n)$ permutes the $i$-th row of the matrix $M(\underline{\mu})=(\mu_{i,j})$ corresponding to $\underline{\mu}$. Thus 
$$
\sum_{\sigma \in R}(-1)^{\sigma}\sigma \cdot\underline{\mu}(\T \cdot c_{n})=\prod^{n}_{i=1}\left[ \sum_{\sigma \in \mathbb{S}_{m}} (-1)^{\sigma}\sigma \cdot\underline{\mu_{i}}(\T \cdot c_{n})\right],
$$
with each $\sum_{\sigma \in \mathbb{S}_{m}} (-1)^{\sigma}\sigma \cdot\underline{\mu_{i}}(\T \cdot c_{n})$ being the Weyl numerator of $\GL(m,\C)$ with $\underline{\mu_{i}}=(\mu_{i,j})_{j}$, evaluated at $(\omega^{i-1}_{n}t_1,\omega^{i-1}_{n}t_2,\cdots,\omega^{i-1}_{n}t_m)$ which is 
\begin{align*}
\prod^{m}_{j=1}(t_{j}\omega^{i-1}_n)^{\mu_{i,j}}&=\left [ \prod^{m}_{j=1}t_j^{\mu_{i,j}}(\omega^{i-1}_n)^{\mu_{i,j}} \right] \\
&=\omega^{m(i-1)^{2}}_{n}\left [ \prod^{m}_{j=1}t_j^{\mu_{i,j}} \right].
\end{align*}
As every $\mu_{i,j}=\mu_{(i-1)m+j} \equiv (i-1) \pmod n$, we write $\mu_{i,j}=n\eta_{i,j}+i-1$, for which 
$$
\prod^{m}_{j=1}(t_{j}\omega^{i-1}_n)^{\mu_{i,j}}=\omega^{m(i-1)^{2}}_{n}\left [\prod^{m}_{j=1}(t^{n}_{j})^{\eta_{i,j}} \right] \left[ \prod^{m}_{j=1}t_{j}\right]^{i-1}.
$$
Therefore 
\begin{align}
\nonumber
\underline{\mu}(\T \cdot c_{n})&=\omega^{\sum^{n}_{i=1}m(i-1)^{2}}_{n}\times \prod^{n}_{i=1}\left [\prod^{m}_{j=1}(t^{n}_j)^{\eta_{i,j}}\right] \times \left [ \prod^{m}_{j=1}t_{j} \right]^{\frac{n(n-1)}{2}} \\
&=\omega^{\sum^{n}_{i=1}m(i-1)^{2}}_{n}\times  \prod^{n}_{i=1}\eta_{i}(\T^{n}) \times \left [ \prod^{m}_{j=1}t_{j} \right]^{\frac{n(n-1)}{2}},\label{eqn 3.7}
\end{align}
where $\eta_{i}=(\eta_{i,j})_{j}$ is the highest weight of $\GL(m,\C)$ defined by:
$$
\underline{\mu_{i}}=n\left (\eta_{i}+\rho_{m}\right)+\underbrace{(i-1,i-1,\cdots,i-1)}_{m \ {\rm times}},
$$
for each $1 \leq i \leq n$. This shows that we have 
\begin{equation} \label{eqn 3.8} 
\sum_{\sigma \in R}(-1)^{\sigma}\sigma\cdot\underline{\mu}(\T \cdot c_{n})=\omega^{\sum^{n}_{i=1}m(i-1)^{2}}_{n}\prod^{n}_{i=1}\left (\sum_{\sigma \in \mathbb{S}_{m}}(-1)^{\sigma}\sigma\cdot\eta_{i}(\T^{n})\right)\left [ \prod^{m}_{j=1}t_{j} \right]^{\frac{n(n-1)}{2}}.
\end{equation}
By Lemma \ref{lem 3} we have 
\begin{align*}
\sum_{\sigma \in R}(-1)^{g \sigma}g \sigma \cdot\underline{\mu}(\T \cdot c_n)&=(-1)^{g}\sum_{\sigma \in R}(-1)^{\sigma}g \sigma\cdot \underline{\mu}(\T \cdot c_n) \\
&=(-1)^{g}\omega^{-m\sum^{n}_{i=1}(i-1)^{2}}_{n}C_{g}(\underline{\mu})\sum_{\sigma \in R}(-1)^{\sigma}\sigma \cdot\underline{\mu}(\T \cdot c_n).
\end{align*}
Now the proof of Proposition \ref{prop 3.6} follows directly from equation \eqref{eqn 3.8}.
\end{proof}
\begin{theorem}\label{4.1}
For an irreducible highest weight representation $\pi_{\underline{\lambda}}$ of ${\rm GL}(mn,\C)$, 
such that for each $i$, $1 \leq i \leq n$, there are exactly $m$ integers in $\underline{\lambda}+\rho_{mn}$ that are congruent to $i-1$ modulo $n$, 
we have 
\begin{equation}\label{Ch Formula} 
\Theta_{\underline{\lambda}}(\T \cdot c_n)=\pm \prod_{i=1}^{n}\Theta_{\eta_i}(\T^{n}), 
\end{equation} 
where $\Theta_{\eta_i}$ is the  character of the highest weight representation of ${\rm GL}(m,\C)$ with highest weight $\eta_{i}$ (see Proposition \ref{prop 3.6}).   
\begin{proof}
From Proposition \ref{prop 3.6} and Lemma \ref{lem 3}, the Weyl numerator is
\begin{align*}
A_{\M}(\T \cdot c_n) &= \sum_{g \in C}\left(\sum_{\sigma \in R}(-1)^{g \sigma}g \sigma \cdot\M(\T \cdot c_n) \right) \\
&= \sum_{g \in C}(-1)^{g}C_{g}(\M) \cdot \left (\prod_{i=1}^{n} S_i(\T^{n})\right) \cdot \left(\prod_{j=1}^{m}t_j\right)^{\frac{n(n-1)}{2}} \\
&= E \cdot \left(\prod^{m}_{j=1}t_{j}\right)^{\frac{n(n-1)}{2}} \cdot \prod^{n}_{i=1}S_{i}(\T^{n}),
\end{align*} 
where $E=\pm \prod_{k<l}(\omega^{k}_{n}-\omega^{l}_{n})^{m}$ and $S_{i}$ are the Weyl numerators for the representations of $\GL(m,\C)$ of highest weight $\eta_{i}$, for $1 \leq i \leq n$.

Further, we have the Weyl denominator
\[
A_{\rho_{mn}}(\T \cdot c_n) = \pm \prod_{k<l}(\omega^{k}_{n}-\omega^{l}_{n})^{m} \cdot \left(\prod^{m}_{j=1}t_{j}\right)^{\frac{n(n-1)}{2}} \cdot \left(\prod_{i<j}(t^{n}_i-t^{n}_{j})\right)^{n}.
\]
Therefore,
\begin{align*}
\frac{A_{\M}(\T \cdot c_n)}{A_{\rho_{mn}}(\T \cdot c_{n})} &= \pm \frac{\prod^{n}_{i=1}S_{i}(\T^{n})}{\left(\prod_{i<j}(t^{n}_i-t^{n}_{j})\right)^{n}} \\
&= \pm \prod^{n}_{i=1}\Theta_{\eta_{i}}(\T^{n}).
\end{align*}
This completes the proof of Theorem \ref{4.1}.
\end{proof}
\end{theorem}

\noindent{\bf Acknowledgment:} The author thanks Prof. Dipendra Prasad for  suggesting the question and for his help in the writing of this paper. The author also thanks the reviewer for valuable comments and helpful suggestions.

\end{document}